\documentclass[11pt]{amsart}

\usepackage{amsmath,amsthm,amsfonts,amssymb,bm,graphicx}

\usepackage{hyperref}
\hypersetup{colorlinks=true,citecolor=blue,linkcolor=blue,urlcolor=blue,
pdfstartview=FitH, pdfauthor=Valentin Blomer and Djordje Mili\'cevi'c ,pdftitle=Kloosterman sums in residue classes}

\theoremstyle{plain}
\newtheorem{theorem}{Theorem}
\newtheorem{cor}[theorem]{Corollary}
\newtheorem{lemma}{Lemma}

\numberwithin{equation}{section}

\theoremstyle{definition}

\renewcommand{\geq}{\geqslant}
\renewcommand{\leq}{\leqslant}

\def\sumstar{\mathop{\sum\nolimits^{\ast}}}

\makeatletter
\DeclareRobustCommand\widecheck[1]{{\mathpalette\@widecheck{#1}}}
\def\@widecheck#1#2{%
    \setbox\z@\hbox{\m@th$#1#2$}%
    \setbox\tw@\hbox{\m@th$#1%
       \widehat{%
          \vrule\@width\z@\@height\ht\z@
          \vrule\@height\z@\@width\wd\z@}$}%
    \dp\tw@-\ht\z@
    \@tempdima\ht\z@ \advance\@tempdima2\ht\tw@ \divide\@tempdima\thr@@
    \setbox\tw@\hbox{%
       \raise\@tempdima\hbox{\scalebox{1}[-1]{\lower\@tempdima\box
\tw@}}}%
    {\ooalign{\box\tw@ \cr \box\z@}}}
\makeatother

\begin{document}

\author{Valentin Blomer}
\author{Djordje Mili\'cevi\'c}
\address{Mathematisches Institut, Bunsenstr.~3-5, D-37073 G\"ottingen, Germany} \email{blomer@uni-math.gwdg.de}
\address{Max-Planck-Institut f\"ur Mathematik, Vivatsgasse 7, D-53111 Bonn, Germany}
\curraddr{Bryn Mawr College, Department of Mathematics, 101 North Merion Avenue, Bryn Mawr, PA 19010, U.S.A.}\email{dmilicevic@brynmawr.edu}

\title{Kloosterman sums in residue classes}


\begin{abstract} We prove upper bounds for sums of Kloosterman sums against general arithmetic weight functions. In particular, we obtain power cancellation in sums of Kloosterman sums over arithmetic progressions, which is of square-root strength in any fixed primitive congruence class up to bounds towards the Ramanujan conjecture.
\end{abstract}


\setcounter{tocdepth}{2}  \maketitle

\section{Introduction} The distribution of values of complete exponential sums is of central interest in number theory and arithmetic geometry. In particular, many arithmetic problems can be transformed into bounding sums of Kloosterman sums. While Weil's bound gives the best possible estimate for individual Kloosterman sums $S(m, n, c)$, one can often use the Bruggeman-Kuznetsov formula to obtain additional savings from the sum over the modulus $c$. Starting with the ground-breaking work of Deshouillers-Iwaniec \cite{DI}, this has been a recurring theme in analytic number theory, see e.g.\ \cite{BFI, DFI1, FM} for some spectacular examples. Using the Bruggeman-Kuznetsov formula for the congruence subgroup $\Gamma_0(q)$, one can require additional divisibility conditions on the modulus $c$. It is a routine exercise to obtain non-trivial bounds for sums of the type
\begin{equation}\label{1}
  \sum_{c \equiv 0 \, (q)} \frac{S(m, n, c)}{c^{1/2}}  f_{\infty}\left(\frac{c}{X}\right),
\end{equation} 
where $f_{\infty} : (0, \infty) \rightarrow \Bbb{C}$ is an appropriate fixed weight function, $m$, $n$ are positive integers, and $X$ is a large parameter.
It is much less of a routine exercise to obtain results of the same quality when the congruence condition $c \equiv 0$ (mod $q$) is replaced with  $c \equiv a$ (mod $q$) for some $(a, q) = 1$. The difficulty  here lies in  the fact that there is no obvious subgroup of $\text{SL}_2(\Bbb{Z})$ where the set of lower left entries of its elements is the set of $c \equiv a$ (mod $q$), and it is also not obvious how to use different cusps of $\Gamma_0(q)$ to encode the congruence condition. It is therefore not clear if spectral theory can provide any non-trivial information even for a fixed progression such as $c \equiv 2$ (mod 5).   In this article we show that nevertheless the problem can be solved in full generality in the framework of certain congruence subgroups $\Gamma_1(Q)$ (more precisely, for some divisors $Q$ of $q^2$).

In fact, we will consider the following more general setup. Let $q$ be an arbitrary positive integer. For a function $f : (\Bbb{Z}/q\Bbb{Z})^{\ast} \rightarrow \Bbb{C}$ we denote by $\widehat{f}$ its ``$q$-Mellin transform"
\begin{displaymath}
  \widehat{f}(\chi) = \frac{1}{\phi(q)^{1/2}} \left.\sum_{c \, (q)}\right.^{\ast} \bar{\chi}(c ) f(c )
\end{displaymath}
where $\chi$ is a Dirichlet character modulo $q$. It satisfies the inversion formula,
\begin{equation}\label{1a}
  f(c ) = \frac{1}{\phi(q)^{1/2}} \sum_{\chi \, (q)} \widehat{f}(\chi) \chi( c).
\end{equation}
We write $\| \widehat{f} \|_1 := \sum_{\chi\,(q)} |\widehat{f}(\chi)|$, and we lift $f$ to a function on integers $c \in \Bbb{Z}$ with $(c, q) = 1$ in the obvious way. We think of $f$ as an \emph{arithmetic weight function} against which we want to sum Kloosterman sums. The original problem of Linnik~\cite{Lin} is concerned with the magnitude of a variant of \eqref{1} with $q=1$ and $f_{\infty}$ a sharp cutoff function; see \cite{Ku}, \cite{ST}. Weighting contributions of various $S(m,n,c)$ by an arithmetic weight $f(c)$ is the natural adelic counterpart of this question, with various moduli $c$ entering with weights according to their position relative to various $p$-adic neighborhoods in addition to the archimedean ones. See \cite{FKM} for a very interesting discussion of arithmetic weights in a different context.\\

Our principal result is the following theorem.

\begin{theorem}
\label{MainTheorem}
Let $m, n, q$ be   positive integers  and $X\geqslant 1$.  Let $f : (\Bbb{Z}/q\Bbb{Z})^{\ast} \rightarrow \Bbb{C}$, and let $f_{\infty} : (0, \infty) \rightarrow \Bbb{C}$ be a   smooth, compactly supported function. Then uniformly in $mn \leq X^2$ one has
\begin{equation}\label{thm1}
  \sum_{(c, q) = 1} \frac{S(m, n, c)}{c^{1/2}}  f(c ) f_{\infty}\left(\frac{c}{X}\right) \ll_{f_{\infty},  \varepsilon}    X^{1/2+2\theta}  \| \widehat{f} \|_1 (mnq)^{\varepsilon}
\end{equation}
for any $\varepsilon > 0$ and any admissible exponent $\theta$ for the Generalized Ramanujan Conjecture for the places dividing $mn$ and the archimedean place.
\end{theorem}

The large range $mn \leq X^2$ of uniformity, also known as the ``Linnik range'' \cite{ST}, is  natural and will become apparent in \eqref{defXi} below. One can also treat the complementary range, but then the analysis of the relevant integral transforms changes. We remark that, if one assumes the Selberg eigenvalue conjecture, then the factor $X^{2\theta}$ in Theorem~\ref{MainTheorem} can be replaced with $(mn)^{\theta}$. Currently, the best available result toward the Generalized Ramanujan Conjecture is $\theta=7/64$, due to Kim and Sarnak~\cite{KS}.

The norm $\|\widehat{f}\|_1$ that appears in Theorem~\ref{MainTheorem} satisfies a general (and generally sharp) estimate $\| \widehat{f} \|_1 \leq \phi(q)^{1/2} \| f \|_2$ by the Cauchy-Schwarz inequality. A particularly interesting class of arithmetic weights $f$ (for $q$ prime) comes from algebraic geometry (e.g.\ as Frobenius trace functions of perverse $l$-adic sheaves \cite{FKM}); see \cite{Ka} for bounds on various norms of the corresponding $\widehat{f}$ in terms  of the ``conductor" of the associated sheaf.

In the proof of Theorem \ref{MainTheorem} we develop a generalized version of the Kuznetsov formula which encodes an additional arithmetic weight  and furnishes an exact spectral decomposition of the left hand side of \eqref{thm1}. We state this spectral formula in Section \ref{sec6} as Theorem \ref{general} after we have developed the necessary notation.  Several important ingredients in the proof of Theorem \ref{MainTheorem} may be noteworthy, including the encoding of arbitrary arithmetic weights using twisted Kloosterman sums in Section~\ref{EncodingSection}, a more general treatment of newforms along the lines of \cite{ILS} in Section~\ref{FourierSection}, and a completely explicit computation of the Fourier coefficients (also at ramified places) of Eisenstein series of general central character in Section \ref{EisensteinSection}.  \\

As a particular application of Theorem~\ref{MainTheorem}, we can choose $f$ to be the characteristic function of the arithmetic progression $c \equiv a$ (mod $q$), where $(a, q) = 1$, and obtain the following variant of  \eqref{1}, a version of Linnik's conjecture \cite{Lin} in arithmetic progressions:

\begin{cor}
\label{coro}
Under the same assumptions as in Theorem~\ref{MainTheorem} one has
\begin{displaymath}
   \sum_{c \equiv a \, (q)} \frac{ S(m, n, c)}{c^{1/2}}  f_{\infty}\left(\frac{c}{X}\right) \ll_{f_{\infty}, \varepsilon}  X^{1/2+2\theta} q^{1/2 } (mnq)^{\varepsilon}.
\end{displaymath}
 \end{cor}

Using Weil's bound $|S(m,n,c)|\leqslant\tau(c)(m,n,c)^{1/2}c^{1/2}$ \cite[Corollary 11.12]{IK} individually, one obtains an upper bound $\ll_{m, n} X^{\varepsilon} (1+X/q)$   in the situation of Corollary \ref{coro}. We see that, for fixed $m, n$, Corollary \ref{coro} is non-trivial for $q\ll X^{(1-4\theta)/3-\varepsilon}$, the quality of cancellation obtained being  uniform across all primitive classes  $a\pmod q$.\\

The dependence on $f_{\infty}$ in Theorem~\ref{MainTheorem} is completely explicit in the proof, see \eqref{mathM} -- \eqref{mathH} below, but it is somewhat convoluted, so that we decided not to display it in the main theorem. Suffice it to say  that Theorem~\ref{MainTheorem} is a very close non-archimedean analogue of the classical situation where $f$ is trivial, but $f_{\infty}$ is oscillating, and our method puts the archimedean and non-archimedean aspects on equal footing. In particular the appearance of norms of $\widehat{f}$ on the right hand side of \eqref{thm1} is rather natural as a comparison with \eqref{mathM} shows.

As an example of this parallelism in Theorem~\ref{MainTheorem}, if we take $f_{\infty}$ fixed and $f(c ) = e(ac/q)$ for some $(a, q) = 1$ and $q$ prime, then $\| \widehat{f} \|_1 = q+\text{O}(1)$. On the other hand, if we take $q=1$, but a similarly oscillating   weight function $f_{\infty}(x) = w(x) e({\tt q}x)$ for some fixed, smooth, compactly supported weight function $w$ and some real number ${\tt q} > 1$, then a stationary phase computation shows that the (usual) Mellin transform of $f_{\infty}$ satisfies $\widehat{f}_{\infty}(c + it) \ll_{ c}  {\tt q}^{-1/2} (1+|t|/{\tt q})^{-10}$, and hence \eqref{mathM} gives the bound  $X^{1/2+2\theta} {\tt q} (mn{\tt q})^{\varepsilon}$. 

An equally robust analogy can be observed in Corollary~\ref{coro}, in which the Kloosterman sums $S(m,n,c)$ are summed over moduli $c$ in a non-archimedean ball \emph{away from $0$}. The archimedean analogue of Corollary \ref{coro} is the situation $q=1$, but a smooth test function $f_{\infty}$ with support in a short interval $[1, 1+{\tt q}^{-1}]$. Then $f_{\infty}$ satisfies $\widehat{f}_{\infty}(c + it) \ll_{ c}  {\tt q}^{-1} (1+|t|/{\tt q}^2)^{-10}$, and \eqref{mathM}  returns the bound $X^{1/2+2\theta} {\tt q}^{1/2} (mn{\tt q})^{\varepsilon}$, in complete analogy with Corollary \ref{coro}. This conclusion should be contrasted with a sum over the progression $c\equiv 0\pmod q$, whose analogue would be the substantially easier case of the support of $f_{\infty}$ in $[0,{\tt q}^{-1}]$ (or, equivalently, having $X/{\tt q}$ in place of $X$). The (perhaps at first counterintuitive) phenomenon that a thinner or shorter sequence gives rise to a harder estimation problem away from the zero class is not uncommon; compare e.g. an application of Voronoi summation to obtain bounds for $\sum_{n \equiv a \, (q)} \lambda_f(n) w(n/X)$ for a fixed cusp form $f$  and $q > X^{1/2}$.

The deeper reason for the  analogy between the archimedean and non-archimedean aspect is 
   that the spectral decomposition of the sum in \eqref{thm1}   using the Bruggeman-Kuznetsov formula gives rise to a spectral sum of Maa{\ss} forms    of a comparable number of terms in both cases:  in slightly simplified terms, we obtain  in the non-archimedean case   a sum with bounded spectral parameter for the group $\Gamma_1(q)$ (of covolume $\asymp q^2$ in Weyl's law), while in the archimedean case we obtain a sum with spectral parameter of size up to ${\tt q}$ (containing $\asymp {\tt q}^2$ terms by Weyl's law)  for $\text{SL}_2(\Bbb{Z})$. In other words, we are expanding in \emph{different  directions} of the full automorphic spectrum of $\text{SL}_2(\mathbb{Q})\setminus\text{SL}_2(\mathbb{A}_{\mathbb{Q}})$. In both cases we estimate the spectral sum trivially and therefore obtain results of comparable quality. (We remark on the side that this numerical phenomenon holds in higher rank, too: one expects $\asymp{\tt q}^{(n^2+n-2)/2}$ Maa{\ss} forms for $\text{SL}_n(\Bbb{Z})$ in a ball of radius ${\tt q}$ about the origin in the Lie-algebra $i\mathfrak{a}^{\ast}$ and $\asymp q^{(n^2+n-2)/2}$ Maa{\ss} forms for $\Gamma_1(q) \subseteq\text{SL}_n(\Bbb{Z})$ in a fixed ball.)\\

We note that, as a direct consequence (see \cite{Va}) of Corollary \ref{coro}, we obtain the following equidistribution result for the Dedekind sums $s(d,c)$. For a real number $x$, let $\langle x\rangle$ denote the fractional part of $x$. 
\begin{cor} 
Let $q$ be a natural number, and let $a$ be an integer coprime to $q$. Then
 the set
$$ \{\langle12 \cdot s(d, c)\rangle  : d\,(\textnormal{mod }c),\,  c \leq  x,\,  c \equiv  a\,  (\textnormal{mod }q) \}$$
becomes equidistributed in $[0, 1)$, as $x \rightarrow  \infty$. 
\end{cor}

Many variations of the present approach are possible. Depending on the application, one can include  additional divisibility conditions on $c$  in \eqref{thm1} and thereby relax the condition $(a, q) = 1$ in Corollary \ref{coro}, one can take $m$ and $n$ to be of opposite sign (using the ``opposite sign" Kuznetsov formula), or, perhaps most interestingly, one 
can sum over $m$ and $n$ and prove large sieve type  inequalities as in \cite{DI}. An investigation of  \eqref{mathM} and in particular the dependence on $f_{\infty}$ in \eqref{thm1} makes it also possible to replace the smooth summation condition by a sharp cutoff condition $c \leq X$. Finally, the dependence on $\theta$ in Theorem \ref{MainTheorem} can be improved slightly in certain ranges of $m, n, q, X$ using density results for exceptional eigenvalues as in \cite[(16.61)]{IK}; see also \cite[(16.75)]{IK}. 
We leave these and other extensions to the interested reader. 

For the rest of paper, implicit constants may depend on $\varepsilon$ (whose numerical value may change from line to line), but all other dependencies are explicitly specified.  \\
 
We would like to thank Farrell Brumley for insightful conversations on topics related to this work and Gergely Harcos for helpful feedback. Henryk Iwaniec kindly informed us that in unpublished notes he obtained cancellation in the situation of Corollary \ref{coro}. Finally we thank the referee for very useful comments. 

The first author acknowledges the support  by the Volkswagen Foundation and a Starting Grant of the European Research Council. The second author would like to thank the Max Planck Institute for Mathematics in Bonn for their support and excellent research infrastructure at the Institute.

\section{Encoding of the arithmetic weights}
\label{EncodingSection}

Let  $\chi_1$ be a primitive Dirichlet character modulo $q_1$, let $m, \delta \in \Bbb{N}$ be positive integers satisfying $(m\delta, q_1) = 1$,  and let $h:\mathbb{N}\to\mathbb{C}$ be a function such that $|h(c)|\ll c^{-3/2-\eta}$ for some $\eta>0$. For every $c$ with $q_1 \mid c$, let
\begin{displaymath}
  S_{\chi_1}(m, n, c) = \left.\sum_{d \, (c )} \right.^{\ast} \chi_1(d) e\left(\frac{md + n \bar{d}}{c}\right)
\end{displaymath}
be the twisted Kloosterman sum, where ${}^{\ast}$ denotes summation restricted to primitive residue classes. By M\"obius inversion, we have that
\begin{equation}
\label{BasicMobius}
   \sum_{(c, q_1) = 1} S_{\chi_1}(m, nq_1^2, q_1\delta c)h(\delta c )  = \sum_{d \mid q_1} \mu(d) \sum_{dq_1 \mid c} S_{\chi_1}(m, nq_1^2, \delta c ) h(\delta c/q_1). 
\end{equation}

On the other hand, if $c$ is such that $(q_1, \delta c) = 1$, we have by twisted multiplicativity of Kloosterman sums (see \cite[(1.59)]{IK}; our case is \emph{mutatis mutandis}) that
\[ S_{\chi_1}(m, nq_1^2, q_1\delta c) = S(m\bar{q}_1, n q_1, \delta c) S_{\chi_1}(m \overline{\delta c}, 0, q_1) =  S(m, n  , \delta c) \bar{\chi}_1(m)\chi_1(\delta c ) \tau(\chi_1) \]
where $\tau(\chi_1)$ is the Gau{\ss} sum, and $\bar{x}$ in one of the first two arguments of a Kloosterman sum denotes the multiplicative inverse of $x$ to the respective modulus. Substituting into \eqref{BasicMobius}, we obtain  
\begin{displaymath}
 \sum_{ c} \chi_1(\delta c ) S(m, n  , \delta c) h(\delta c )=  \frac{\chi_1(m)}{\tau(\chi_1)}\sum_{d \mid q_1} \mu(d) \sum_{dq_1 \mid c} S_{\chi_1}(m, nq_1^2, \delta c ) h(\delta c/q_1).
\end{displaymath}

Now, let $\chi$ be an arbitrary Dirichlet character modulo $q$, induced by a primitive character $\chi_1$ modulo $q_1 \mid q$, let $(m,q_1)=1$, and let $h$ be as above. Then
\begin{displaymath}
\begin{split}
   \sum_{ (c, q) = 1} \chi (  c ) S(m, n  ,  c) h(  c )& = \sum_{\substack{\delta \mid q\\ (\delta, q_1) = 1}} \mu(\delta) \sum_{c} \chi_1 (  \delta c ) S(m, n  ,  \delta c) h( \delta  c )\\
    & = \sum_{\substack{\delta \mid q\\ (\delta, q_1) = 1}} \mu(\delta) \frac{\chi_1(m)}{\tau(\chi_1)}\sum_{d \mid q_1} \mu(d) \sum_{\delta dq_1 \mid c} S_{\chi_1}(m, nq_1^2,   c ) h(  c/q_1).
  \end{split} 
\end{displaymath}
 We can collapse the double sum over $d$ and $\delta$ to a single sum getting an equality
\begin{equation}\label{basic}
   \sum_{ (c, q) = 1} \chi (  c ) S(m, n  ,  c) h(  c )  =    \frac{\chi_1(m)}{\tau(\chi_1)}\sum_{d \mid q} \mu(d) \sum_{ dq_1 \mid c} S_{\chi_1}(m, nq_1^2,   c ) h(  c/q_1),
\end{equation}
valid for every $m$ with $(m,q_1)=1$ (and so a fortiori whenever $(m,q)=1$).

For general $m$, we put $m' = m/(m, q^{\infty})$ and $n' = n(m, q^{\infty})$. Since $S(m, n, c) = S(m', n', c)$ for $(c, q) = 1$, \eqref{basic} holds without the condition $(m, q) = 1$, if we replace $m$ and $n$ by $m'$ and $n'$ on the right hand side; note that $\chi(m')=\chi_1(m')$. We thus obtain the more general equality
\begin{equation}\label{basicgeneral}
\sum_{(c,q)=1}\chi(c)S(m,n,c)h(c)= \frac{\chi(m')}{\tau(\chi_1)}\sum_{d\mid q}\mu(d)\sum_{ dq_1\mid c}S_{\chi_1}(m',n'q_1^2,c)h(c/q_1),
\end{equation}
valid for any positive integers $m, n, c$ and any Dirichlet character $\chi$ modulo $q$. Here, note that $mn = m'n'$. 

Equality \eqref{basicgeneral} is at the heart of our argument. We think of the left-hand side of this formula as an average of Kloosterman sums in the usual sense (as a smooth sum over $c$) but additionally weighted with a special arithmetic weight, namely an arbitrary multiplicative character of the modulus $c$. Our equality expresses such a twisted average in terms of twisted Kloosterman sums, which we \emph{can} analyze using spectral theory of automorphic forms for $\Gamma_0(dq_1)$ and character $\chi_1$, i.e.\ a character of the quotient $ \Gamma_1(dq_1)\backslash \Gamma_0(dq_1)$. In this sense, as remarked in the introduction, we solve the problem of bounding Kloosterman sums with arithmetic weights modulo $q$ (with arbitrary weights in the next paragraph) in the framework of suitable congruence subgroups  $\Gamma_1(dq_1)$, where $d,q_1 \mid q$. Recall that the additional factor $d$ was inherited simply from M\"obius inversion; it will turn out to be of little relevance in the forthcoming asymptotic analysis.

Finally, we encode an arbitrary arithmetic weight $f:(\mathbb{Z}/q\mathbb{Z})^{\ast}\to\mathbb{C}$ and a long-range archimedean cutoff. With $f$ and $f_{\infty}$ as in the statement of Theorem~\ref{MainTheorem}, we use \eqref{basicgeneral} with $h(c)=f_{\infty}(c/X)/c^{1/2}$, multiply by $\widehat{f}(\chi)$, sum over all Dirichlet characters $\chi$ modulo $q$, and use the inversion formula \eqref{1a}. Thus  we obtain the basic identity
\begin{equation}
\label{THEbasicidentity}
\begin{split}
&\sum_{(c, q) = 1}\frac{S(m, n, c) }{c^{1/2}}f(c ) f_{\infty}\left(\frac{c}{X}\right)\\
&=\frac{1}{\phi(q)^{1/2}} \sum_{\chi \, (q)} \frac{\chi(m') \widehat{f}(\chi) }{\tau(\chi_1)}   \sum_{d \mid q} \mu(d) \sum_{dq_1 \mid c} \frac{q_1^{1/2}}{c^{1/2}} S_{\chi_1}(m', n'q_1^2, c ) f_{\infty}\left(\frac{c}{q_1X}\right). 
\end{split}
\end{equation}
Our identity relates the sum of Kloosterman sums against both a finite and an archimedean test function (essentially an almost arbitrary smooth, compactly supported function on $\mathbb{Q}^{\times}\setminus\mathbb{A}_{\mathbb{Q}}^{\times}$) to sums which may be treated by the Kuznetsov trace formula for the group $\Gamma_0(dq_1)$. Here we can also see the underlying motivation for the sums appearing in \eqref{BasicMobius}.

This prepares ground for our principal application, Theorem~\ref{MainTheorem}. In the following sections, we will prove that
\begin{equation}\label{claim}
\Sigma_{\chi_1}(m,n,d,q,X):=\sum_{dq_1 \mid c} \frac{1}{c^{1/2}} S_{\chi_1}(m', n'q_1^2, c ) f_{\infty}\left(\frac{c}{q_1X}\right) \ll_{f_{\infty}}  q_1^{1/2} X^{1/2+2\theta} (mnq)^{\varepsilon},
\end{equation}
uniformly in $mn \leq X^2$ and across all $d,q_1\mid q$ and all primitive characters $\chi_1$ modulo $q_1$. Taking this for granted, Theorem~\ref{MainTheorem} follows from trivial estimates.  We note that, if one is not interested in  uniformity with respect to $m, n, q$, then  \cite[Theorem 1]{GS} shows directly that
$$\Sigma_{\chi_1}(m,n,d,q,X) \ll_{m, n, q, f_{\infty} } X^{1/2+2\theta+\varepsilon}.$$ 
Remark 1 in \cite{GS} gives some explicit polynomial  dependence on $m, n, q$, but considerably weaker than  required for \eqref{claim}. If  $\chi_1$ is trivial, one can also read off uniform bounds for \eqref{claim} from  \cite[(16.72), (16.75)]{IK}. Modulo the Ramanujan-Petersson conjecture, these bounds have an extra factor $(n(m, q^{\infty})q_1)^{1/4}$ compared to \eqref{claim}. In the following section we will make systematic use of newform theory in order to optimize the dependence on $n, m$ and $q_1$. 
 Our general treatment of newforms and Eisenstein series and corresponding bounds, in particular \eqref{ortho} and Lemma~\ref{EisensteinLemma} below, as well as Theorem~\ref{general}, may be of independent interest.
 
\section{Fourier coefficients of automorphic forms}
\label{FourierSection}

In this section, which can be read independently of Section~\ref{EncodingSection}, we collect facts and conventions about (holomorphic and Maa{\ss}) cusp forms and Eisenstein series and present estimates for their Fourier coefficients which will be used in the estimation of $\Sigma_{\chi_1}$ in Section~\ref{KuznetsovSection}.

Let $\chi_1$ be a primitive character modulo $q_1$, let $\kappa = 0$ if $\chi_1$ is even and $\kappa = 1$ if $\chi_1$ is odd, and let $d$ be a square-free integer. (In our application to \eqref{claim}, $q_1$ and $d$ will have the same meaning as in the rest of the paper.) 
For a positive integer $k \geq 2$ satisfying $k\equiv\kappa\pmod 2$, let $\mathcal{A}_k(dq_1,\chi_1)=S_k(\Gamma_0(dq_1),\chi_1)$ denote the finite-dimensional space of holomorphic weight $k$ cusp forms of level $dq_1$ and character $\chi_1$. Here and on, for any $r$ with $q_1\mid r\mid dq_1$, $\chi_1$ in the notation for a space such as $\mathcal{A}_k(r,\chi_1)$ (or a basis such as $\mathcal{B}_k$ below) stands, strictly speaking, for the induced character $\chi_1\chi_0$, where $\chi_0$ is the principal character modulo $r$. For simplicity, we suppress $\chi_0$ from notation (but not from computation) as no confusion will arise. Let $\mathcal{A}_{\kappa}(dq_1,\chi_1)=L^2_{\text{cusp}}(\Gamma_0(dq_1)\setminus G,  \chi_1, \kappa)$ denote the space of non-holomorphic weight $\kappa$ cusp forms on $G=\text{SL}_2(\mathbb{R})$ of level $dq_1$ and character $\chi_1$, and, for every $t \in(\mathbb{R} \cup [-i/2, i/2])/\{\pm 1\}$, let $\mathcal{A}_{\kappa}(dq_1,\chi_1, t)\subset\mathcal{A}_{\kappa}(dq_1,\chi_1)$
denote the finite-dimensional space of forms in $\mathcal{A}_{\kappa}(dq_1,\chi_1)$ of Laplacian eigenvalue $1/4 + t^2$.  

We normalize Hecke operators  both for holomorphic and non-holomorphic cusp forms so that the Ramanujan-Petersson conjecture states that the eigenvalues $\lambda_f(p )$ of $T_p$ are bounded by 2 in absolute value (i.e.\ \cite[(6.1)]{DFI} for Maa{\ss} forms and   \cite[(2.15)]{ILS}   for holomorphic forms of   weight $k$).

With respect to the standard inner product 
\begin{equation}
\label{inner}
\langle f, g\rangle
= \int_{\Gamma_0(dq_1)\backslash \mathfrak{h}} f(z) \bar{g}(z)y^{\ell} \frac{dx\, dy}{y^2}\end{equation}
of level $dq_1$, where $\ell=k$ in the space $\mathcal{A}_k(dq_1,\chi_1)$ and $\ell=0$ in each of the spaces $\mathcal{A}_{\kappa}(dq_1,\chi_1)$, we construct specific orthonormal bases $\mathcal{B}_{k}(dq_1, \chi_1)$ and $\mathcal{B}_{\kappa}(dq_1, \chi_1)$ of these spaces as follows.  For every $r$ satisfying $q_1\mid r \mid dq_1$, let $\mathcal{A}_k^{\ast}(r,\chi_1)$ denote the space of  holomorphic  forms orthogonal to all oldforms of level $r$ and character $\chi_1$. By Atkin-Lehner theory  (in particular the multiplicity one principle), we can choose an orthonormal basis $\mathcal{B}_k^{\ast}(r, \chi_1)$  of $\mathcal{A}_k^{\ast}(r,\chi_1)$ consisting of \emph{newforms}, i.e. eigenforms for \emph{all} Hecke operators $T_n$ ($n \geq 1$)  with eigenvalue $\lambda_f(n)$, say   (see \cite[Section 14.7]{IK}).  For $f \in \mathcal{B}_k^{\ast}(r, \chi_1)$, let $\mathcal{A}_{dq_1/r}(f)$ be the space spanned by the set of shifts $\mathcal{S}_{dq_1/r}(f):=\{f (bz) : b \mid dq_1/r\}$. Then we have (again by multiplicity one) an orthogonal decomposition
\[ \mathcal{A}_k(r,\chi_1)=\bigoplus_{q_1\mid r\mid dq_1}\bigoplus_{f\in\mathcal{B}_k^{\ast}(r,\chi_1)}\mathcal{A}_{dq_1/r}(f). \]
One can obtain an orthonormal basis of $\mathcal{A}_{dq_1/r}(f)$ by orthonormalizing the set $\mathcal{S}_{dq_1/r}(f)$; below, we will
specify an explicit orthonormal basis $\mathcal{S}_{dq_1/r}^{\square}(f)=\{f_{(b)}(z) : b \mid dq_1/r\}$. We obtain the requisite orthonormal basis $\mathcal{B}_k(dq_1,\chi_1)$ of the entire space $\mathcal{A}_k(dq_1,\chi_1)$ by taking the union of all these sets:
\[ \mathcal{B}_k(dq_1,\chi_1):=\bigsqcup_{q_1\mid r\mid dq_1}\bigsqcup_{f\in\mathcal{B}_k^{\ast}(r,\chi_1)}\mathcal{S}_{dq_1/r}^{\square}(f). \]
We construct an orthonormal basis $\mathcal{B}_{\kappa}(dq_1, \chi_1, t)$ of the space $\mathcal{A}_{\kappa}(dq_1,\chi_1, t)$ analogously and write
\begin{displaymath}
 \mathcal{B}_{\kappa}(dq_1, \chi_1) = \bigsqcup_t \mathcal{B}_{\kappa}(dq_1, \chi_1, t),
\end{displaymath}
the union being taken over the   spectral resolution of $\mathcal{A}_{\kappa}(dq_1,\chi_1)$. 

We write the Fourier expansion of a modular form $f$ as
\begin{equation}
\label{fourier}
\begin{alignedat}{2}
f(z)&=\sum_{n\geqslant 1}\rho_f(n)n^{k/2}e(nz)&\quad&\text{for }f \in \mathcal{A}_k(dq_1, \chi_1),\\
f(z) &=\sum_{n\neq 0}\rho_f(n)W_{\frac n{|n|}\frac{\kappa}2,it_f}\big(4\pi|n|y\big)e(nx)&&\text{for }f \in \mathcal{A}_{\kappa}(dq_1, \chi_1, t_f),
\end{alignedat}
\end{equation}
 where $W_{\alpha, \beta}$ is the Whittaker function. Let $f$ be a newform of level $r$ and character $\chi_1$. The Fourier coefficients are related to the Hecke eigenvalues as 
\begin{equation}\label{hecke}
  \sqrt{n} \rho_f(n) = \rho_f(1) \lambda_f(n)
\end{equation}  
for all $n \in \Bbb{N}$. Moreover, if $f$ is  $L^2$-normalized  with respect to \eqref{inner}, then we have the following bounds that are essentially  due to Hoffstein and Lockhart (upper bound) and Iwaniec (lower bound):
\begin{equation}\label{rho1}
|\rho_f(1)|^2 = \begin{cases}  \displaystyle \frac{\cosh(\pi t_f)}{  dq_1 (1+ |t_f|)^{\kappa}}  (dq_1(1+|t_f|))^{o(1)}, &  f\in\mathcal{A}_{\kappa}^{\ast}(r,\chi_1, t_f),\\[5mm]
  \displaystyle\frac{(4\pi)^{k-1}}{dq_1\Gamma(k)}(kdq_1)^{ o(1)}, & f\in\mathcal{A}_k^{\ast}(r,\chi_1);\end{cases}
\end{equation}
see \cite[(30) -- (31)]{HM}, \cite{HL}. (Compare the slightly different normalization in \cite[Section 2.2]{HM}.)  We recall the Hecke relations
\begin{equation}
\label{heckerelations}
\begin{alignedat}{2}
&\lambda_f(p^{\alpha+1}) = \lambda_f( p)\lambda_f(p^{\alpha}) - \chi_1(p)\chi_0(p)\lambda_f(p^{\alpha-1})&\quad&\text{for all primes $p$ and integers $\alpha \geq 1$,}\\
&\chi_1( p) \bar{\lambda}_f(p ) = \lambda_f( p)&& \text{for primes }p \nmid r,
\end{alignedat}
\end{equation}
(see e.g. \cite[p.\ 371]{IK} or \cite[p.\ 520]{DFI}),  as well as the bounds
\begin{equation}
\label{bounds}
\begin{alignedat}{2}
&  |\lambda_f(p )| \leq p^{\theta} + p^{-\theta}, &\quad& p \nmid r,\\
&  |\lambda_f(p  )| \leq 1,&& p \mid r;
\end{alignedat}
\end{equation}
cf.\  \cite[Theorem 1.1 ii) and iii)] {Li} (as well as our normalization \eqref{fourier} and  \eqref{hecke}) for the latter bound, which holds verbatim for Maa{\ss} forms. 

For our principal estimation in the next section, we want to be completely explicit in our construction of the basis $\mathcal{S}_{dq_1/r}^{\square}(f)$ and compute the Fourier coefficients of the various $f_{(b)}$. For notational simplicity we consider only  the holomorphic case; the Maa{\ss} case is identical upon replacing $k$ with $0$ throughout the argument. We follow closely the  argument in \cite[Section 2]{ILS} (see also \cite[Section 3]{AU}), which requires only minor modification.   Let $f \in \mathcal{B}_k^{\ast}(r, \chi_1)$ and recall that $dq_1/r$ (and thus each of its divisors $b$) is squarefree. Define the arithmetic function
\begin{displaymath}
  \nu(b) := b  \prod_{p \mid b} \left(1+ \frac{\chi_0(p )}{p}\right)
\end{displaymath}
where $\chi_0$ is the principal character modulo $r$, 
and write $f|_{b}(z) := b^{k/2} f(bz)$.
 Starting from the expression  
\begin{displaymath}
 \langle E(\cdot, s) f(b_1 \cdot {}), f(b_2 \cdot {})\rangle =  \int_{\Gamma_0(dq_1)\backslash \mathfrak{h}} E(z, s) f(b_1z) \bar{f}(b_2z) y^k \frac{dx\, dy}{y^2}
 \end{displaymath}
 where $E(z, s)$ is the standard weight 0 non-holomorphic Eisenstein series of level $dq_1$, unfolding, using \eqref{heckerelations} to explicitly evaluate multiplicatively shifted convolution $L$-series as a product of Euler factors
\begin{displaymath}
\sum_{\alpha=0}^{\infty} \lambda_f(p^{\alpha+1}) \bar{\lambda}_f(p^{\alpha}) p^{-\alpha s} =  \lambda_f( p)(1 + \chi_0( p)p^{-s})^{-1}\sum_{\alpha=0}^{\infty} |\lambda_f( p^{\alpha})|^2 p^{-\alpha s},
\end{displaymath} and evaluating residues at $s=1$, we find as in \cite[Lemma 2.4]{ILS} that
\begin{displaymath}
  \langle f|_{b_1}, f|_{b_2}\rangle = \frac{\bar{\lambda}_f(b')\lambda_f(b'')}{\nu(b') \nu(b'')} (b'b'')^{1/2} \langle f, f\rangle, \quad \quad b' = \frac{b_1}{(b_1, b_2)}, \quad b'' = \frac{b_2}{(b_1, b_2)}. 
\end{displaymath}
The Maa{\ss} case is identical except that slightly different special functions occur in the unfolding step; see \cite[Section 19]{DFI}. 
With this at hand, proceeding as in \cite[p.\ 75]{ILS}, we find that the forms
\begin{equation}\label{gram}
  f_{(b)}(z) = \Bigl\{b \prod_{p \mid b} \Bigl( 1 - \frac{p|\lambda_f(p )|^2}{(p + \chi_0(p ))^2}\Bigr)^{-1}\Bigr\}^{1/2} \sum_{c \ell = b} \frac{\mu(c )\bar{\lambda}_f(c )}{\nu(c )} \ell^{(k-1)/2} f(\ell z)
\end{equation}
form an orthonormal basis of $\mathcal{A}_{dq_1/r}(f)$; this is the basis $\mathcal{S}_{dq_1/r}^{\square}(f)$ of our choice. (Taking inverses in the product above is justified by \eqref{bounds} with any $\theta < 1/2$.) 
 Hence for every $f\in\mathcal{B}^{\ast}_k(r,\chi_1)$ or $f\in\mathcal{B}^{\ast}_{\kappa}(r,\chi_1, t)$, every $n \in \Bbb{N}$ and every $b\mid(dq_1/r)$, and with $f_{(b)}$ as in \eqref{gram}, we have that 
\begin{equation}\label{ortho}
\begin{split}
\sqrt{n}  \rho_{f_{(b)}}(n)& =  \Bigl\{b \prod_{p \mid b} \Bigl( 1 - \frac{p|\lambda_f(p )|^2}{(p + \chi_0(p ))^2}\Bigr)^{-1}\Bigr\}^{1/2} \sum_{c\ell = b} \frac{\mu(c )\bar{\lambda}_f(c )}{\nu(c )} \rho_f(1)\lambda_f\left(\frac{n}{\ell}\right)\\
& \ll b^{\varepsilon} \sum_{\ell \mid b} \frac{\ell \,|\lambda_f(b/\ell)|}{b^{1/2}} \left|\rho_f(1)\lambda_f\left(\frac{n}{\ell}\right)\right|,
\end{split}
\end{equation}
with the convention that $\lambda_f(x) = 0$ for $x \not\in \Bbb{N}$. Here, we only used so far that \eqref{bounds} holds with any $\theta<1/2$, and $|\rho_f(1)|$ can be further estimated by \eqref{rho1}.
\\

We proceed to a discussion of the Eisenstein spectrum, which is parametrized by singular cusps $\mathfrak{a}$. Write $Q=dq_1$. Recall that a cusp $\mathfrak{a}$ for a group $\Gamma$ is called singular with respect to a multiplier system $\vartheta$ on $\Gamma$ if $\vartheta(\gamma)=1$ for all $\gamma$ in the stabilizer $\Gamma_{\mathfrak{a}}\subset\Gamma$. For a cusp $\mathfrak{a}$ for the group $\Gamma=\Gamma_0(Q)$, denoting by $\sigma_{\mathfrak{a}}$ a scaling matrix (that is, a matrix such that $\sigma_{\mathfrak{a}}^{-1}\Gamma_{\mathfrak{a}}\sigma_{\mathfrak{a}}=\left\{\pm\left(\begin{smallmatrix} 1&k\\0&1\end{smallmatrix}\right):k\in\mathbb{Z}\right\}$), $\mathfrak{a}$ is a singular cusp for a character $\chi$ modulo a divisor of $Q$ if\footnote{See \cite[p.\ 44]{IwClassical}; compare with \cite[(4.43)]{DFI}. It is tempting to think of the two elements $\pm\sigma_{\mathfrak{a}}\left(\begin{smallmatrix} 1&1\\&1\end{smallmatrix}\right)\sigma^{-1}_{\mathfrak{a}}$ as playing the same role, but this is not quite so, due to the presence of the factor $j_{\tau}(z)^{-\kappa}$ in the multiplier.}
\[ \chi\left(\sigma_{\mathfrak{a}} \left(\begin{matrix} 1 & 1\\ & 1\end{matrix}\right) \sigma_{\mathfrak{a}}^{-1}\right) =1. \]
As usual, we interpret $\chi$ as a character on $\Gamma_0(Q)$ via $\chi(\gamma) =  \chi(d) = \bar{\chi}(a)$ for $\gamma = \left(\begin{matrix} a & b\\ c & d\end{matrix}\right)$.
 
For a singular cusp $\mathfrak{a}$, we consider the weight $\kappa$ Eisenstein series
 \begin{equation}\label{defE}
  E_{\mathfrak{a}, \chi_1} (z, s) = \sum_{\gamma \in \Gamma_{\mathfrak{a}}\backslash \Gamma} \bar{\chi}_1(\gamma) j_{\sigma_{\mathfrak{a}}^{-1}\gamma}(z)^{-\kappa} \Im(\sigma_{\mathfrak{a}}^{-1}\gamma z)^{s} = \sum_{\tau \in \Gamma_{\infty} \backslash \sigma_{\mathfrak{a}}^{-1}\Gamma} \bar{\chi}_1(\sigma_{\mathfrak{a}}\tau) j_{\tau}(z)^{-\kappa}\Im (\tau z)^{s},
\end{equation}
where $j_{\tau}(z) = \frac{\tilde{c}z+\tilde{d}}{|\tilde{c}z+\tilde{d}|}$ for $\tau =\left(\begin{smallmatrix} \tilde{a} & \tilde{b}\\ \tilde{c} & \tilde{d}\end{smallmatrix}\right)$ is the usual multiplier, and $\tau=\sigma_{\mathfrak{a}}^{-1}\gamma$. This series converges absolutely for $\Re s > 1$, and it has a meromorphic extension to all of $\Bbb{C}$. Eisenstein series have a Fourier expansion similar to Maa{\ss} forms, which at the point $1/2 + it$ is given by
\begin{equation}
\label{FourierDefE}
   E_{\mathfrak{a}, \chi_1} (z, 1/2 + it) = C_{\mathfrak{a},\chi_1, t}(z)+ \sum_{n\neq 0}\rho_{\mathfrak{a}, \chi_1}(n, t)W_{\frac n{|n|}\frac{\kappa}2,it}\big(4\pi|n|y\big)e(nx).
\end{equation}

In Section~\ref{EisensteinSection}, we specify a full list of inequivalent singular cusps $\mathfrak{a}$ in \eqref{gcd}, explicitly compute the Fourier coefficients $\rho_{\mathfrak{a},\chi_1}(n, t)$ of the Eisenstein series $E_{\mathfrak{a},\chi_1}(z, 1/2 + it)$ in \eqref{double}, and prove the following uniform upper bound, which will be used in our principal estimation in the next section.

\begin{lemma}
\label{EisensteinLemma}
Let $m,n,Q$ be positive integers, let $\chi_1$ be a character modulo $Q$, let $t$ be a real number, and let $\kappa\in\{0,1\}$. Let $\tilde{Q}$ be the smallest positive integer such that $Q\mid\tilde{Q}^2$. Then, for every cusp $\mathfrak{a}$ of $\Gamma=\Gamma_0(Q)$ singular for the character $\chi_1$, the Fourier coefficients $\rho_{\mathfrak{a},\chi_1}(n,t)$ in the expansion \eqref{FourierDefE} of the weight $\kappa$ Eisenstein series $E_{\mathfrak{a},\chi_1}(z,s)$ defined by \eqref{defE} satisfy
\begin{displaymath}
\mathcal{E}(m,n,t):=\sum_{\mathfrak{a} \text{ {\rm singular}}}\frac{ \sqrt{mn}}{\cosh(\pi t)} \overline{\rho_{\mathfrak{a}, \chi_1}(m, t) } \rho_{\mathfrak{a}, \chi_1}(n, t) \ll  \frac{1}{(1+|t|)^{\kappa}}\frac{(m,\tilde{Q})^{1/2}(n,\tilde{Q})^{1/2}}{\tilde{Q}}(Qmn(1+|t|))^{\varepsilon}.
\end{displaymath}
\end{lemma}
We note that the second fraction does not exceed 1 and that $\tilde{Q} \geq Q^{1/2}$. We postpone this more detailed discussion of the Eisenstein series and the proof of Lemma~\ref{EisensteinLemma} to the end of the paper. We will not concern ourselves with the constant term $C_{\mathfrak{a},\chi_1, t}(z)$ since it does not enter our estimates. 

\section{Application of the Kuznetsov formula}
\label{KuznetsovSection}

In this section, we use the Kuznetsov trace formula and the estimates from Section~\ref{FourierSection}, including \eqref{ortho} and Lemma~\ref{EisensteinLemma}, to prove our claim \eqref{claim}. With notation from Section~\ref{EncodingSection}, put
\begin{equation}\label{defg}
  g(x) := \Bigl(\frac{4\pi \sqrt{m'n'q_1^2}}{x}\Bigr)^{1/2} f_{\infty}\Bigl( \frac{4\pi \sqrt{m'n'}}{x X}\Bigr).
\end{equation}  
Then $g(x)$ is a smooth function compactly supported on an interval of $x$ satisfying \begin{equation}\label{defXi}
  x \asymp \Xi := \sqrt{mn} X^{-1} \leq 1
 \end{equation}
   and such that $\| g \|_{\infty} \asymp (q_1X)^{1/2}\|f\|_{\infty}$, and we have $$c^{1/2} f_{\infty}\Bigl(\frac{c}{q_1X}\Bigr) = g\Bigl(\frac{4\pi \sqrt{m'n'q_1^2}}{c}\Bigr).$$
Substituting into the left-hand side of \eqref{claim}, we have that
\[ \Sigma_{\chi_1}(m,n,d,q,X)=\sum_{dq_1\mid c}\frac1cS_{\chi_1}(m',n'q_1^2,c)g\Bigl(\frac{4\pi\sqrt{m'n'q_1^2}}c\Bigr). \]
The sum on the right-hand side can now be readily transformed with the Kuznetsov formula for the group $\Gamma_0(dq_1)$ and character $\chi_1$, which we quote from Blomer-Harcos-Michel~\cite{BHM}. We use the usual weight 0 or the weight 1 formula according as $\chi_1$ is even or odd. 
We define the following integral transforms:
\begin{displaymath}
\begin{split}
  \dot{g}(k) & = i^k \int_0^{\infty} J_{k-1}(x) g(x) \frac{dx}{x}, \quad k \in \Bbb{N},\\
  \tilde{g}(t) & = \frac{it^{\kappa}}{2\sinh(\pi t)} \int_0^{\infty} \bigl(J_{2it}(x) - (-1)^{\kappa} J_{-2i t}(x)\bigr) g(x) \frac{dx}{x}, \quad t \in \Bbb{R} \cup [-i/2, i/2].
\end{split}
\end{displaymath}
The power series expansion of the Bessel functions \cite[8.440]{GR} together with \eqref{defXi} yields
\begin{equation}\label{gdot}
  \dot{g}(k) \ll \| g \|_{\infty} \Gamma(k)^{-1}
\end{equation}
and 
\begin{equation}\label{gtilde}
  \tilde{g}(t)\ll  \begin{cases}  |t|^{\kappa} \left(\displaystyle\frac{|\widehat{g}(2 i t)| + |\widehat{g}(-2 i t)|}{|t|^{1/2}} + \frac{|\widehat{g}(2+2 i t)| + |\widehat{g}(2-2 i t)|}{|t|^{3/2}}  + \frac{\| g \|_{\infty}}{|t|^{5/2}} \right), & |t|\geqslant 1,\\[0.3cm]
     \| g \|_{\infty} \Xi^{-2\theta}, & t \in (-1, 1) \cup [-i\theta, i \theta]. \end{cases}
\end{equation}\\
Here $\widehat{g}$ denotes the Mellin transform of $g$, and by \eqref{defg} we have
\begin{equation}\label{mellin}
  \widehat{g}(s) = (q_1X)^{1/2} \left(\frac{4\pi \sqrt{mn}}{X}\right)^s \widehat{f}_{\infty}(1/2 - s). 
\end{equation}

The following spectral summation formula holds \cite[p.\ 705]{BHM}:
\begin{equation}
\label{SpectralSummation}
     \sum_{dq_1 \mid c} \frac{1}{c} S_{\chi_1}(m', n'q_1^2, c )g\Bigl(\frac{4\pi \sqrt{m'n'q_1^2}}{c}\Bigr) = \mathcal{H} + \mathcal{M} + \mathcal{E},
\end{equation}   
   where
   \begin{displaymath}\begin{split}
& \mathcal{H} =   \underset{\substack{k \equiv \kappa\, (2), k > \kappa\\ f \in \mathcal{B}_{k}(dq_1, \chi_1)}}{\sum\sum} \dot{g}(k) \frac{(k-1)! \sqrt{m'n'q_1^2}}{\pi (4\pi)^{k-1}} \overline{\rho_f(m')} \rho_f(n'q_1^2), \\
 & \mathcal{M}=   \sum_{f \in \mathcal{B}_{\kappa}(dq_1, \chi_1)} \tilde{g}(t_f) \frac{4\pi \sqrt{m'n'q_1^2}}{\cosh(\pi t_f)}  \overline{\rho_f(m')} \rho_f(n'q_1^2), \\
& \mathcal{E} =  \sum_{ \mathfrak{a} \text{ singular}} \int_{-\infty}^{\infty} \tilde{g}(t) \frac{ \sqrt{m'n'q_1^2}}{\cosh(\pi t)}  \overline{\rho_{\mathfrak{a}, \chi_1}(m', t)} \rho_{\mathfrak{a}, \chi_1}(n'q_1^2, t)\, dt.
  \end{split} 
\end{displaymath}
In comparison with \cite{BHM}, note that we are using the classical parametrization of the Eisenstein spectrum in terms of singular cusps. The formula \eqref{SpectralSummation} is purely spectral in that $\mathcal{B}_k(dq_1,\chi_1)$ and $\mathcal{B}_{\kappa}(dq_1,\chi_1)$ can be arbitrary orthonormal bases of the spaces $\mathcal{A}_k(dq_1,\chi_1)$ and $\mathcal{A}_{\kappa}(dq_1,\chi_1)$, respectively. However, we will from now on assume that these bases have been chosen as in Section~\ref{FourierSection}, which will allow us to efficiently estimate the terms on the right-hand side.

We start with bounding $\mathcal{M}$, the contribution of the Maa{\ss} spectrum. Each summand in $\mathcal{M}$ corresponds to a basis vector of the form $f_{(b)}$ for some newform $f\in\mathcal{B}_{\kappa}^{\ast}(r,\chi_1)$, with $q_1\mid r\mid dq_1$ and $b\mid dq_1/r$.

Since $(m', q) = (m', dq_1) = 1$,  \eqref{bounds} and \eqref{ortho} imply
\[ \sqrt{m'} \rho_{f_{(b)}}(m') \ll b^{-1/2+\varepsilon} |\rho_f(1) \lambda_f(m' b)|\ll (m')^{\theta+\varepsilon}b^{\theta-1/2+\varepsilon}|\rho_f(1)|. \]
On the other hand, writing $b_0=(b,r^{\infty})$, $b_1=b/b_0$, $n'_0=(n'q_1^2,r^{\infty})$, and $n'_1=n'q_1^2/n'_0$, we have that
\begin{align*}
\sqrt{n'q_1^2}\rho_{f_{(b)}}(n'q_1^2)
&\ll b^{-1/2+\varepsilon}|\rho_f(1)|\sum_{\ell_0\mid(b_0,n'_0)}\sum_{\ell_1\mid (b_1,n'_1)}\ell_0\ell_1\left(\frac{b_1}{\ell_1}\right)^{\theta+\varepsilon}\left(\frac{n'_1}{\ell_1}\right)^{\theta+\varepsilon}\\
&\ll (n'_1)^{\theta+\varepsilon}b^{-1/2+\varepsilon}b_0b_1^{1-\theta}|\rho_f(1)|\ll (n')^{\theta+\varepsilon}b^{1/2+\varepsilon}|\rho_f(1)|.
\end{align*}
Combining these estimates with \eqref{rho1}, we  obtain the following bound for any individual term occuring in the sum for $\mathcal{M}$:
\begin{equation}\label{maass}
  \frac{4\pi \sqrt{m'n'q_1^2}}{\cosh(\pi t_f)}  \overline{\rho_{f_{(b)}}(m')} \rho_{f_{(b)}}(n'q_1^2) \ll \frac{(mnb)^{\theta}}{dq_1 (1+|t_f|)^{\kappa-\varepsilon}}
  (qmn)^{\varepsilon} \ll  \frac{(mn)^{\theta}}{r (1+|t_f|)^{\kappa-\varepsilon}}
  (qmn)^{\varepsilon}.
\end{equation}
(Note that we can afford to let go of a factor of $b^{1-\theta}$.)
Next we  use the well-known fact that
\begin{equation}\label{weyl}
  \#\{ f \in \mathcal{B}^{\ast}_{\kappa}(r, \chi_1) :  |t_f| \leq T \} \ll
  (r T^2)^{1+\varepsilon}
\end{equation}
for any $r$ with $q_1 \mid r$. This weak but very uniform version of Weyl's law follows for instance by combining \eqref{rho1} with \cite[(16.56)]{IK} in the special case $n=1$. 

Collecting \eqref{defXi},  \eqref{gtilde}, \eqref{mellin}, \eqref{maass} and  \eqref{weyl}, denoting by $\tau(\cdot)$ the divisor function, and using summation by parts, we conclude that
\begin{equation}\label{mathM}
\begin{split}
   \mathcal{M} &\ll q_1^{1/2} X^{1/2+2\theta} (mnq)^{\varepsilon}
   \sum_{q_1\mid r\mid dq_1}\frac{\tau(dq_1/r)}r
   \Bigg[ \#\big\{f\in\mathcal{B}^{\ast}_{\kappa}(r,\chi_1):|t_f|<1\big\}\cdot\|f_{\infty}\|_{\infty}+\\
&\mskip 240mu   \sum_{\substack{f\in\mathcal{B}^{\ast}_{\kappa}(r,\chi_1)\\ |t_f|\geqslant 1}} \Bigl(\frac{|\widehat{f}_{\infty}(1/2 \pm 2 i t_f)|}{|t_f|^{1/2-\varepsilon}} +\frac{|\widehat{f}_{\infty}(-3/2 \pm 2 i t_f)|}{|t_f|^{3/2-\varepsilon}}+\frac{\|f_{\infty}\|_{\infty}}{|t_f|^{5/2-\varepsilon}}\Bigr)\Bigg]\\
   & \ll_{f_{\infty}}  q_1^{1/2} X^{1/2+2\theta} (mnq)^{\varepsilon},
 \end{split}  
 \end{equation}   
as required for \eqref{claim}.  

For the Eisenstein spectrum we replace \eqref{maass} with the bound given in Lemma~\ref{EisensteinLemma} for $Q = dq_1$. Recalling that $(m', dq_1) = 1$,  we obtain a slightly stronger estimate 
\begin{equation}\label{mathE}
\begin{split}
 \mathcal{E}& =\int_{-\infty}^{\infty}\tilde{g}(t)\mathcal{E}(m',n'q_1^2,t)dt\ll \frac{q_1^{1/2} X^{1/2 }}{(dq_1)^{1/4}} (mnq)^{\varepsilon} \bigg(\|{f_{\infty}}\|_{\infty} + \int_{|t| \geq 1} \frac{|\widehat{f}_{\infty}(1/2 - 2 i t)|}{|t|^{1/2-\varepsilon}} dt \bigg)\\
 &  \ll_{f_{\infty}} \frac{q_1^{1/2} X^{1/2 }}{(dq_1)^{1/4}} (mnq)^{\varepsilon} . 
 \end{split}
\end{equation} 
The bound for the holomorphic spectrum $\mathcal{H}$ is along the same lines with \eqref{gdot} instead of \eqref{gtilde}, giving
\begin{equation}\label{mathH} \mathcal{H} \ll \|{f_{\infty}}\|_{\infty}\cdot q_1^{1/2} X^{1/2 } (mnq)^{\varepsilon}. \end{equation}
This completes the proof of Theorem~\ref{MainTheorem}.
 
\section{Proof of Lemma~\ref{EisensteinLemma}}
\label{EisensteinSection}

It remains to prove Lemma~\ref{EisensteinLemma}. To this end we compute  explicitly  the Fourier coefficients $\sqrt{n} \rho_{\mathfrak{a}, \chi_1}(n, t)$. This is in principle straightforward, but a bit tedious.

We follow \cite{DI}, which treats the case when $\chi_1$ is trivial.   We first describe a set of (inequivalent) \emph{singular} cusps. Given a divisor $w\mid Q$ and a primitive residue class $r$ modulo $w_Q:=(w,Q/w)$, we can always find a residue class $u\bmod w$ such that $u\equiv r\pmod{w_Q}$ and $(u,w)=1$ (since the coprimality condition $(u,p)=1$ is automatic from $u\equiv r\pmod{w_Q}$ for primes $p\mid w_Q$ and can be imposed by the Chinese Remainder Theorem by requiring, for example, that $u\equiv 1\pmod p$ for $p\nmid w_Q$). Let $\mathcal{U}_w$ be a full set of representatives $u$ of primitive residue classes modulo $w_Q$ chosen so that $(u,w)=1$ for every $u\in\mathcal{U}_w$. Then, a full set of inequivalent cusps for the group $\Gamma = \Gamma_0(Q)$ is given by all fractions 
  $\mathfrak{a} = \frac{u}{w}$, where $w  \mid Q$ and $u\in\mathcal{U}_w$; see \cite[Lemma 2.3]{DI}. 
A possible scaling matrix for a cusp $\mathfrak{a}$ is given by (see \cite[p.\ 247]{DI})
\begin{displaymath}
  \sigma_{\mathfrak{a}} = \left(\begin{matrix} u \sqrt{Q/(w^2, Q)} & 0 \\ w \sqrt{Q/(w^2, Q)} & \frac{1}{u\sqrt{Q/(w^2, Q)}}\end{matrix}\right).
\end{displaymath}
 We compute that, for every $k\in\mathbb{Z}$,
\begin{displaymath}
  \sigma_{\mathfrak{a}} \left(\begin{matrix} 1 & k\\ & 1\end{matrix}\right) \sigma_{\mathfrak{a}}^{-1} = \left(\begin{matrix}1 - kuwQ/(w^2, Q) & ku^2Q/(w^2, Q)\\ - kw^2Q/(w^2, Q) & 1 + kuwQ/(w^2, Q)  \end{matrix}\right).
\end{displaymath}
Hence the singular cusps $\mathfrak{a}=u/w$ for the character $\chi_1$ of conductor $q_1$ are given by the condition 
\begin{equation}\label{gcd}
  (w, Q/w) \mid \frac{Q}{q_1}  \quad \Longleftrightarrow \quad q_1 \mid [w, Q/w].
\end{equation}
This description of the set of equivalence classes of singular cusps can also be found in \cite[Lemma 13.5]{IwClassical}. 

Recall the definition \eqref{defE} of the Eisenstein series at a singular cusp $\mathfrak{a}=u/w$ as above. A convenient parametrization of the sum on the right-hand side of \eqref{defE} is as follows. For $\gamma = \left(\begin{smallmatrix} a & b\\ c & d\end{smallmatrix}\right) \in \Gamma$, we have
\begin{align*}
\sigma_{\mathfrak{a}}^{-1}\gamma=\sigma_{\mathfrak{a}}^{-1}  \left(\begin{matrix} a & b\\ c & d\end{matrix}\right)
&= \left(\begin{matrix} \frac{a}{u \sqrt{Q/(w^2, Q)}} & \frac{b}{u \sqrt{Q/(w^2, Q)}}\\ (uc - wa) \sqrt{Q/(w^2, Q)} & (ud-wb)\sqrt{Q/(w^2 , Q)}\end{matrix} \right)\\
&=:   \left(\begin{matrix} \ast  & \ast\\ -Cw \sqrt{Q/(w^2, Q)} & D\sqrt{Q/(w^2 , Q)}\end{matrix} \right),
\end{align*} 
where $Cw=-uc+wa$ and $D=ud-wb$ satisfy $\gamma^{-1} \cdot \frac{u}{w} = \frac{D}{Cw}$. In particular, as $\gamma$ runs through $\Gamma_{\mathfrak{a}}\setminus\Gamma$, the point $D/Cw$ precisely traverses the orbit of the cusp $u/w$ in $\Gamma$. By \cite[Lemma 3.6]{DI} the set of pairs $(Cw,D)$ is characterized by 
\begin{displaymath}
  (Cw, D) = 1, \quad (C, Q/w) = 1, \quad  CD \equiv u \, \text{mod } (w, Q/w).
\end{displaymath}
Pairs of integers $(C,D)$ with these properties come in couples $\pm(C,D)$. To each such pair $(C,D)$ with $C>0$ thus corresponds a unique class $\gamma\in\Gamma_{\mathfrak{a}}\setminus\Gamma$; its representatives $\left(\begin{smallmatrix} a&b\\c&d\end{smallmatrix}\right)$ clearly satisfy $a \equiv C \, (\text{mod }Q/w)$ and $d \equiv \bar{u}D \, (\text{mod } w)$.

By \eqref{gcd}, the character $\chi_1$ can be induced from the product $\psi_1\psi_2$, where $\psi_1$ is a primitive character of some conductor $r_1 \mid w$ and $\psi_2$ is a primitive character of some conductor $r_2 \mid Q/w$; in particular, $\chi_1(d) = \psi_1(\bar{u}D)\psi_2(\bar{C})$.
Proceeding as in \cite[p.\ 247]{DI}\footnote{There is a small typo in \cite[p.\ 247]{DI}: the congruence condition in the right-most sum in the first display under Lemma 3.6 should be $\delta \gamma \equiv u \, ((w, q/w))$. Also note the additional factor $\sqrt{n}$ in \cite[Section 2.1.3]{BHM} due to the different weight function compared to \cite[(1.17)]{DI}.}  or \cite[p.\ 526]{DFI} we conclude that for $n > 0$ the $n$-th Fourier coefficient of $E_{\mathfrak{a}, \chi_1}(z, 1/2 + it)$ is given by
\begin{equation}\label{coeff}
\sqrt{n}\rho_{\mathfrak{a}, \chi_1}(n, t)=\frac{i^{\kappa}\pi^s n^{s-\frac{1}{2}}}{\Gamma(s+\frac{\kappa}{2})}\Bigl( \frac{(w, Q/w)}{wQ}\Bigr)^s \sum_{(C, Q/w) = 1} \frac{\psi_2(C )}{C^{2s}}\underset{\substack{D \, (Cw)\\ CD \equiv u \, ((w, Q/w))}}{\left.\sum \right.^{\ast}} \psi_1(u\bar{D}) e\left(-\frac{n D}{Cw}\right),
\end{equation}
where  $s = 1/2 + it$. 

We transform this expression further into a form which will be convenient for our purposes. We detect the congruence condition in the innermost sum by Dirichlet characters $\rho$ modulo $w_Q$, getting
\begin{displaymath}
  \underset{\substack{D \, (Cw)\\ CD \equiv u \, (w_Q)}}{\left.\sum \right.^{\ast}} \psi_1(u\bar{D}) e\left(-\frac{n D}{Cw}\right) = \frac{\psi_1(u)}{\phi(w_Q)} \sum_{\rho\,(w_Q)} \rho(\bar{u}C)  \underset{ D \, (Cw) }{\left.\sum \right.^{\ast}} \bar{\psi}_1 \rho(D)e\left(-\frac{n D}{Cw}\right).
\end{displaymath}
Let $Q_w$ be the product of all prime factors $p\mid Q$ such that $p\nmid (Q/w)$. The condition that $(C,Q/w)=1$ is equivalent to the statement that $C = C_1C_2$, where $(C_1,Q)=1$ and $C_2\mid Q_w^{\infty}$ (with the correspondence being that $C_2=(C,Q)$ and $C_1=C/C_2$). Recalling that the conductor of $\bar{\psi}_1\rho$ is a divisor of $w$, the innermost sum over $D$ above equals
\[ \sumstar_{d_1\,(C_2w)}\sumstar_{d_2\,(C_1)}\bar{\psi}_1\rho(C_1d_1+C_2wd_2)e\left(-\frac{n(C_1d_1+C_2wd_2)}{C_1C_2w}\right)=r_{C_1}(n)\bar{\psi}_1\rho(C_1)S_{\bar{\psi}_1\rho}(-n,0,C_2w), \]
where $\displaystyle r_{C_1}(n)=\sumstar_{d\,(C_1)}e(-nd/C_1)$ is the Ramanujan sum. Substituting, the double sum over $C$ and $D$ in \eqref{coeff} equals
\[ \frac{\psi_1(u)}{\phi(w_Q)} \sum_{\rho\,(w_Q)} \rho(\bar{u})    \sum_{(C_1,Q)= 1}  \frac{\rho\psi_2(C_1)}{C_1^{2s}} \sum_{C_2\mid Q_w^{\infty}} \frac{\rho\psi_2(C_2)}{C_2^{2s}} r_{C_1}(n) \bar{\psi}_1\rho(C_1) S_{\bar{\psi}_1\rho}(-n, 0, C_2w). \]
The inside sum over $C_1$ equals
\[ \sum_{(C_1,Q)=1}\frac{\bar{\psi}_1\psi_2\rho^2(C_1)}{C_1^{2s}}\sum_{\delta\mid(C_1,n)}\mu\left(\frac{C_1}{\delta}\right)\delta=\sum_{\substack{\delta\mid n\\ (\delta,Q)=1}}\frac{\bar{\psi}_1\psi_2\rho^2(\delta)}{\delta^{2s-1}}\frac1{L^{(Q)}(2s,\bar{\psi}_1\psi_2\rho^2)}, \]
where $L^{(Q)}(s,\psi)=\prod_{p\nmid Q}\big(1-\psi(p)p^{-s}\big)^{-1}$ is the partial $L$-function. Putting everything together in \eqref{coeff}, we have that
\begin{equation}\label{double}
\begin{split}
\sqrt{n}\rho_{\mathfrak{a},\chi_1}(n,t)&=\frac{i^{\kappa} \pi^{\frac{1}{2}+it} n^{it}}{\Gamma(\frac{1+\kappa}2+it)}\left(\frac{w_Q}{wQ}\right)^{\frac12+it}\frac{\psi_1(u)}{\phi(w_Q)}\\
&\quad\times\sum_{\rho\,(w_Q)} \frac{\rho(\bar{u})}{L^{(Q)}(1+2it, \bar{\psi}_1\psi_2\rho^2)}  \sum_{\substack{\delta \mid n\\(\delta,Q)=1}} \frac{\bar{\psi}_1\psi_2\rho^2(\delta)}{\delta^{2it}} \sum_{C_2\mid Q_w^{\infty}} \frac{\rho\psi_2(C_2)}{C_2^{1+2it}}   S_{\bar{\psi}_1\rho}(-n, 0, C_2w).
\end{split}
\end{equation}

With this explicit computation of the Fourier coefficients $\sqrt{n}\rho_{\mathfrak{a},\chi_1}(n,t)$ under our belt, we are now ready to prove Lemma~\ref{EisensteinLemma}.   Indeed, we have that
\begin{displaymath}
\mathcal{E}(m,n,t) = \sum_{\substack{w  \mid Q\\ q_1 \mid [w, Q/w]}} \sum_{u \in \mathcal{U}_w} \frac{\sqrt{mn}}{\cosh(\pi t)}  \overline{\rho_{u/w, \chi_1}(m, t) } \rho_{u/w, \chi_1}(n, t).
\end{displaymath}
We insert \eqref{double}, sum over $u$ first, and then estimate trivially, using Stirling's formula, standard lower bounds for $L^{(Q)}(1+it,\psi)$ (including Siegel's bound if $\psi$ is real and $|t| \leq 1$), and Weil's bound for Kloosterman sums. Denoting by $\tau(x)$ the number of divisors of $x$, we have for $m,n> 0$ that
\begin{displaymath} 
\begin{split}
\mathcal{E}(m,n,t)
&\ll\frac{(mn)^{\varepsilon}}{\cosh(\pi t)|\Gamma(\frac{1+\kappa}2+it)|^2}  \sum_{\substack{w  \mid Q\\ q_1 \mid [w, Q/w]}}\frac{w_Q}{wQ\phi(w_Q)}\\
&\qquad\qquad\times \sum_{\rho\,(w_Q)}\frac1{|L^{(Q)}(1+2it,\bar{\psi}_1\psi_2\rho^2)|^2}\mathop{\sum\sum}_{C_2,C_2'\mid Q_w^{\infty}}\frac{|S_{\bar{\psi}_1\rho}(-m,0,C_2w)||S_{\bar{\psi}_1\rho}(-n,0,C'_2w)|}{C_2C'_2}\\
&\ll\frac{(Qmn(1+|t|))^{\varepsilon}}{(1+|t|)^{\kappa}}  \sum_{w \mid Q} \frac{w_Q}{Q}\mathop{\sum\sum}_{C_2,C'_2\mid Q_w^{\infty}}\frac{(m,C_2w)^{1/2}(n,C'_2w)^{1/2}\tau(C_2w)\tau(C'_2w)}{C_2^{1/2}C'_2{}^{1/2}}\\
&\ll\frac{(Qmn(1+|t|))^{\varepsilon}}{(1+|t|)^{\kappa}}  \sum_{w \mid Q}\frac{w_Q(m,w)^{1/2}(n,w)^{1/2}}{Q}\ll  \frac{1}{(1+|t|)^{\kappa}}\frac{(m,\tilde{Q})^{1/2}(n,\tilde{Q})^{1/2}}{\tilde{Q}}(Qmn(1+|t|))^{\varepsilon},
\end{split}  
\end{displaymath} 
as desired.

\section{A generalized Kuznetsov formula}\label{sec6}

Combining \eqref{THEbasicidentity}, \eqref{defg} and \eqref{SpectralSummation}, we obtain the following new version of the Kuznetsov formula. 

\begin{theorem}
\label{general}
Let $m, n, q \in \Bbb{N}$, $f : (\Bbb{Z}/q\Bbb{Z})^{\ast} \rightarrow \Bbb{C}$ and $f_{\infty} : (0, \infty) \rightarrow \Bbb{C}$ be a smooth, compactly supported function. We keep the notation developed so far. In particular,    for a character $\chi$ modulo $q$ we write $\kappa = 0$ if $\chi$ is even and $\kappa = 1$ if $\chi$ is odd, and we denote by $\chi_1$ modulo $q_1$ the underlying primitive character. We write $m' := m/(m, q^{\infty})$ and $n' = n(m, q^{\infty})$. Then
\begin{displaymath}
\begin{split}
  \sum_{(c,q)=1} &\frac{S(m, n, c)}{c^{1/2}} f(c ) f_{\infty}(c )\\
  & = \sum_{d \mid q} \mu(d) \sum_{\chi \, (q)} \sum_{f\in \mathcal{B}_{\kappa}(dq_1, \chi_1)} \frac{4\pi \sqrt{m'n'q_1^2}}{\cosh(\pi t_f)} \chi(m')\overline{\rho_f(m')}\rho_f(n'q_1^2) F(\chi) F_{\infty}(t_f)\\
   & +  \sum_{d \mid q} \mu(d) \sum_{\chi \, (q)}  \sum_{ \substack{\mathfrak{a} \text{ {\rm singular}}\\ \text{\rm{level} } dq_1}} \int_{-\infty}^{\infty}  \frac{ \sqrt{m'n'q_1^2}}{\cosh(\pi t)}  \chi(m')\overline{\rho_{\mathfrak{a}, \chi_1}(m', t)} \rho_{\mathfrak{a}, \chi_1}(n'q_1^2, t)  F(\chi) F_{\infty}(t)\, dt\\
   & +  \sum_{d \mid q} \mu(d) \sum_{\chi \, (q)}  \underset{\substack{k \equiv \kappa\, (2), k > \kappa\\ f \in \mathcal{B}_{k}(dq_1, \chi_1)}}{\sum\sum}  \frac{(k-1)! \sqrt{m'n'q_1^2}}{\pi (4\pi)^{k-1}} \chi(m')\overline{\rho_f(m')} \rho_f(n'q_1^2)  F(\chi) F_{\infty}^{\ast}(k)
\end{split}
\end{displaymath}
where
\begin{displaymath}
\begin{split}
 F(\chi)&  = \frac{\widehat{f}(\chi)q_1^{1/2}}{\tau(\chi_1)\phi(q)^{1/2}},\\
 F_{\infty}(t) & = \frac{\pi^{1/2} i (mn)^{1/4} q_1^{1/2} t^{\kappa} }{\sinh(\pi t)} \int_0^{\infty} (J_{2it}(x) -(-1)^{\kappa}  J_{-2it}(x))   f_{\infty} \left(\frac{4\pi \sqrt{mn}}{x}\right) \frac{dx}{x^{3/2}},\\
  F_{\infty}^{\ast}(k) & = (4\pi)^{1/2} i^k (mn)^{1/4} q_1^{1/2} \int_0^{\infty}  J_{k-1}(x) f_{\infty} \left(\frac{4\pi \sqrt{mn}}{x}\right)  \frac{dx}{x^{3/2}}.
\end{split} 
\end{displaymath}
Here $\mathcal{B}_{\kappa}(dq_1, \chi_1)$ and $ \mathcal{B}_{k}(dq_1, \chi_1)$ can be any orthonormal bases of the spaces of non-holomorphic weight $\kappa$ (holomorphic weight $k$, respectively) cusp forms of level $dq_1$ and character $\chi_1$. 
\end{theorem}

We remark that, despite its appearance, the right hand side of our generalized Kuznetsov formula is symmetric in $m$ and $n$. This is especially easy to see for bases $\mathcal{B}_{\kappa/k}(dq_1,\chi_1)$ consisting of forms $f$ which are eigenforms of Hecke operators $T_m$ with $(m,q)=1$, so that $\sqrt{m}\rho_f(mn)=\lambda_f(m)\rho_f(n)$ whenever $(m,nq)=1$; by \eqref{ortho}, such is the case, in particular, for the special bases constructed in Section~\ref{FourierSection}. In this case, referring also to \eqref{heckerelations}, we have that
\[ \chi(m') \overline{\rho_f(m')}\rho_f(n'q_1^2)  =\sqrt{(mn,q^{\infty})/mn}\lambda_f(m/(m, q^{\infty}))\lambda_f(n/(n, q^{\infty})) \overline{\rho_f(1)}\rho_f((mn, q^{\infty})q_1^2). \]
Our Eisenstein series are in general \emph{not} Hecke eigenfunctions, but the statement of Theorem~\ref{general} (being a purely spectral formula) holds also for the basis of Eisenstein series parametrized by pairs of characters as described in \cite{BHM}; this basis satisfies the usual Hecke  relations, and   the symmetry in the Eisenstein term can be restored by the same argument. \\

\end{document}